\documentclass[english,12p]{article}

    \newcommand{\mytitle}
    {The Weyl group of type $A_1$ root systems extended by an abelian
      group}                                                          

    \usepackage[all]{xy}
    \usepackage{amssymb,amsmath,latexsym}
    \usepackage{eucal}
     \usepackage[
       pdftitle={\mytitle},
       ]{hyperref}
    \usepackage{theorem}

    \newtheorem{lem}{Lemma}[section]
    \newtheorem{thm}[lem]{Theorem}
    
    \newtheorem{cor}[lem]{Corollary}
    {\theorembodyfont{\rm}
      \newtheorem{defn*}[lem]{Definition}
      \newtheorem{rem*}[lem]{Remark}
      \newtheorem{exmp}[lem]{Example}

      \newtheorem{constr*}[lem]{Construction}
    }

    \newenvironment{defn}{\begin{defn*}}{\hspace*{\fill}
        \raisebox{0mm}{$\diamond$}
      \end{defn*}}
    {\par\noindent~\begin{minipage}[b]{0.97\linewidth}}%
    {\end{minipage}}
    \newenvironment{lastalign*}%
    {\par\noindent\begin{minipage}[b]{0.95\linewidth}
        \medskip\begin{center}\begin{math}\displaystyle}%
    {\end{math}\end{center}
          \end{minipage}}
    {\par\medskip\noindent\begin{minipage}[b]{0.95\linewidth}
      \begin{enumerate}}%
    {\end{enumerate}\end{minipage}}
    \newenvironment{proof}{\par\medskip\noindent\upshape\rmfamily
      \textbf{Proof}.\ }{\hspace*{\fill}\rule{1.2ex}{1.2ex}\par\medskip}
    \newenvironment{proof*}{\par\medskip\noindent\upshape\rmfamily
      \textbf{Proof}.\ }{\par\medskip}
    \newenvironment{lastequation*}{
      \par\medskip\noindent\hspace*{\fill}\begin{math}\displaystyle
      }{
      \end{math}\hspace*{\fill}\rule{1.2ex}{1.2ex}\par}
    \newenvironment{lasteqnarray*}{
      \par
      \begin{minipage}[t]{\hfill}
        \begin{eqnarray*}
        }{      
        \end{eqnarray*}
      \end{minipage}
      \begin{minipage}[b]{0.5cm}
        \rule{1.2ex}{1.2ex}
      \end{minipage}
    }
    
    \newcommand{\newop}[2]{\newcommand{#1}{\mathop{\mathsf{\strut #2}}\nolimits}}
    \newop{\Hom}{\mathrm{Hom}}
    \newop{\Aut}{\mathrm{Aut}}
    \newop{\End}{\mathrm{End}}
    \newop{\co}{\mathrm{co}}
    \newop{\id}{\mathrm{id}}
    \newop{\im}{\mathrm{im}}
    \newop{\rk}{\mathrm{rk}}
    \newop{\spann}{\mathrm{span}}
    \newop{\conv}{\mathrm{conv}}
    \newop{\rank}{\mathrm{rank}}
    \newop{\GL}{\mathrm{GL}}
    \newop{\PDS}{\mathrm{PDS}}
    \newop{\Oo}{\mathrm{O}}
    \newop{\cone}{\mathrm{cone}}
    \newop{\algint}{\mathrm{algint}}
    \newop{\inn}{\mathrm{int}}
    \newop{\aff}{\mathrm{aff}}

    \newcommand{\F}{\mathbb{F}}

    \newop{\M}{\mathbf{M}}

    \newop{\PGL}{\mathrm{PGL}}
    \newop{\sL}{\mathrm{sl}}
    \newop{\tr}{\mathrm{tr}}
    \newcommand{\cd}{{\cdot}}
    \newcommand{\sbullet}{\hspace{0.5mm}\begin{picture}(1,5)(-.5,-3)\circle*{2}\end{picture}\hspace{0.5mm}}

    \newcommand{\propname}[1]{{\rm\bf(#1)}}
    \renewcommand{\phi}{\varphi}
    
    \newcommand{\sym}{\mathrm{sym}}

    \newcommand{\ab}{\mathrm{ab}}

    \newcommand{\calU}{\mathcal{U}}
    
    \newcommand{\calV}{\mathcal{V}}

    \newcommand{\calW}{\mathcal{W}}
    
    \newcommand{\calA}{\mathcal{A}}

    \newcommand{\calX}{\mathcal{X}}
    
    \newcommand{\calY}{\mathcal{Y}}

\begin{document}

    \title{\mytitle\footnote{ 
    The research for this article was supported by a postdoctoral
    fellowship at 
    the Department of Mathematics and Statistics at
    Dalhousie University.
    }}
    \author{Georg W. Hofmann}
    
    \maketitle
    \begin{abstract}
      We investigate the class of root systems $R$ obtained by extending an
      $A_1$-type 
      irreducible root system by a free abelian group $G$. 
      In this context there are two reflection groups with respect to
      a discrete symmetric space $T$ associated to $R$, namely, the
      Weyl group $\calW$ of $R$ and a group
      $\calU$ with a so-called presentation by
      conjugation. 
      We show that the
      natural homomorphism $\calU\to\calW$ is an isomorphism if and
      only if an associated subset $T^\ab\setminus\{0\}$ of $G_2=G/2G$ is
      2-independent, i.e. its image under the map
      $G_2\to G_2\otimes G_2,~g\mapsto g\otimes g$ is linearly independent 
      over the Galois field $\F_2$.
      \hfill\break
      \hfill\break                      
    \end{abstract} 
    
    \pagestyle{headings}

    \section{Introduction}                                           %

    We consider 
    a root system $R$ extended by an abelian group $G$, a notion that
    is introduced in \cite{yoshii-rootabelian}. 
    It generalizes the concepts of 
    extended affine root systems (see \cite{EALA_AMS}, for instance) 
    and affine root systems in the sense of
    \cite{saito1}, both of which are generalizations of 
    root systems of affine Kac-Moody
    algebras (see \cite{MoP}, for instance). The Weyl group $\calW$
    of $R$ is not necessarily a Coxeter group, so a more general
    presentation is needed to capture the algebraic structure of $\calW$.
    The group $\calU$ is given by the so-called
    presentation by conjugation with respect to $R$:
    \begin{align*}
      \calU
      \cong
      \big\langle(\hat r_\alpha)_{\alpha\in R^\times}
      ~|~
      &
      \hat r_\alpha=\hat r_\beta
      \text{~~~if $\alpha$ and $\beta$ are linearly dependent,}\\
      &\hat r_\alpha^2=1,~~
      \hat r_\alpha\hat r_\beta\hat r_\alpha^{-1}
      =\hat r_{r_\alpha(\beta)}; 
      \text{ for }\alpha,\beta\in R\big\rangle.
    \end{align*} 
    There is a natural group 
    homomorphism from $\calU$ onto $\calW$. 
    \par\medskip\noindent
    {\bf Question}~
     Is $\calU\to\calW$ injective? In other words, does $\calW$ have
     the presentation by conjugation with respect to $R$?
    \par\medskip\noindent
    This question has been studied for various root systems in
    \cite{krylyuk},
    \cite{EAWG},
    \cite{azamReduced},
    \cite{MR2341017},
    \cite{azamA1},
    \cite{myWeyl},
    \cite{myAbel}. 
    
    In this note we investigate the case that $R$ is of
    type $A_1$, i.e. the underlying finite root system consists of two
    roots. This type of root system $R$ allows for less rigidity then
    other types and is therefor of special
    interest as a prototype. 
    We prove the following result that allows an answer to
    the question above using an algorithmic approach.

    Suppose $R$ is a type $A_1$ 
    root system extended by a free abelian group
    $G$. Then a subset $T^\ab$ of $G_2=G/2G$ can be associated to it
    in a natural way. This subset is called 2-independent, if its
    image under $G_2\to G_2\otimes G_2, g\mapsto g\otimes g$ is a
    linearly independent set. 
    \par\medskip\noindent
    {\bf Theorem}~
     The natural homomorphism 
     \begin{math}
       \calU\to\calW
     \end{math}
     is an isomorphism if and only if $T^\ab\setminus\{0\}$ is
     2-independent in $G_2$.
    \par\medskip\noindent

    This result provides an attractive alternative to a
    characterization proved in \cite{azamA1} using so-called
    integral collections. Our answer to the question
    above is more general than that in \cite {azamA1} 
    as $G$ is not required to
    be finitely generated.
    We expect that the idea of 2-independence that we have introduced
    will play an important role in understanding the question 
    for root systems of the types $B_n$ and $C_n$.

    \section[Discrete symmetric spaces]
    {Discrete symmetric spaces and their reflection groups}  %

    In this section we provide the basic terminology for the following
    sections. The notion of a discrete symmetric space is a special
    case of the symmetric spaces introduced in \cite{MR0239005}. The
    associated category of reflection groups is introduced in
    \cite{myAbel} and more details can be found there.

    \begin{defn}\propname{Discrete symmetric space}
      Let $T$ be a set with a (not necessarily associative) multiplication
      \begin{align*}
        \mu: T\times T\to T,~
        (s,t)\mapsto s\cd t.
      \end{align*}
      Then the pair $(T,\mu)$ is called a
      \emph{discrete symmetric space}\index{discrete symmetric space}
      if the following conditions are satisfied for all $s$, $t$ and
      $r\in T$:
      \begin{enumerate}
  \item[(S1)]
        \begin{math}
          s\cd s=s,
        \end{math}
  \item[(S2)]
        \begin{math}
          s\cd(s\cd t)=t,
        \end{math}
  \item[(S3)]
        \begin{math}
          r\cd(s\cd t)=(r\cd s)\cd(r\cd t).
        \end{math}
      \end{enumerate}
      By abuse of language, 
      we will sometimes say that $T$ is a discrete symmetric space
      instead of saying that $(T,\mu)$ is a discrete symmetric space. 
      If
      $s\cd t=t$ for all $s$ and $t\in T$ then we call $\mu$ the 
      \emph{trivial multiplication}. 
    \end{defn}
    

    For the remainder of this section, let $T$ be a discrete symmetric space. 

    \begin{defn}\propname{Reflection group}\label{defn:reflGroup}
      Let  $\calX$ be a group acting on $T$. 
      We will denote the element in $T$ obtained by $x$ acting on
      $t$ by $x.t$. 
      Let 
      \begin{equation*}
        \sbullet^\calX:~
        T\to\calX,~t\mapsto t^\calX
      \end{equation*}
      be a function. Then $(\calX,\sbullet^\calX)$ is called a
      $T$-\emph{reflection group}%
      \index{reflection group},
      if the following conditions are satisfied:
      \begin{enumerate}
  \item[(G1)] The group $\calX$ is generated by the set
        $T^\calX:=\{t^\calX~|~t\in T\}$.
  \item[(G2)] For all $s$ and $t\in T$ we have
        \begin{math}
          t^\calX.s=t\cd s.
        \end{math}
  \item[(G3)] For all $s$ and $t\in T$ we have
        \begin{math}
          t^\calX
          s^\calX
          (t^\calX)^{-1}
          ~=~
          (t\cd s)^\calX.
        \end{math}
  \item[(G4)] For every $t\in T$ we have
        \begin{math}
          (t^\calX)^2
          =
          1.
        \end{math}
      \end{enumerate}
      If we do not need to specify the map $\sbullet^\calX$ we will
      also say that $\calX$ is a reflection group instead of
      saying that $(\calX,\sbullet^\calX)$ is a reflection
      group.
    \end{defn}

    \begin{defn}\propname{Reflection morphism}
      Let 
      $\calX$ and 
      $\calY$ be $T$-reflection
      groups.
      Then a group
      homomorphism $\phi:\calX\to\calY$ is called a 
      \emph{$T$-morphism},
      if $\phi(t^\calX)=t^\calY$ for every $t\in T$.
    \end{defn}


    Let the group $\calU$ be given by the presentation
    \begin{align*}      
      \calU:=
      \big\langle(t^\calU)_{t\in T}
      ~\big|~
      (t^\calU)^2=1        
      \text{~~and~~}
      t^\calU s^\calU(t^\calU)^{-1}
      =(t\cd s)^\calU
      \text{~~~for $s$ and $t\in T$}
      \big\rangle.
    \end{align*}
    There is map
    $\sbullet^\calU:~T\to\calU,~t\mapsto t^\calU$ associated to the
    presentation. An action of $\calU$ on $T$ can be defined
    satisfying $t^\calU.s=t\cd s$ for all $s$ and $t\in T$. With
    respect to this action the pair $(\calU,\sbullet^\calU)$ is a
    $T$-reflection group.
    There is a unique
    $T$-morphism from $\calU$ into any other $T$-reflection group.

    \begin{defn}
      The pair $(\calU,\sbullet^\calU)$ is called 
      \emph{the initial $T$-reflection group.}      
    \end{defn}
    
    \section[Type $A_1$ root systems extended by an abelian group]
    {Type $A_1$ root systems extended by an abelian group}   %

    In this section we introduce the concept of a type $A_1$ root
    system extended by an abelian group $G$ in an ad hoc manner. Thus
    we avoid presenting the details of the definition for more general
    types.

    Let $(G,+)$ be an abelian group. Define the
    multiplication 
    \begin{align}\label{eq:mult}
      G\times G\to G,~
      (g,h)\mapsto g\cd h=2g-h.
    \end{align}
    Now let $T$ be a generating subset of $G$ such that $0\in T$ and
    $G\cd T\subseteq T$. 
    It is straightforward to verify that $T$ with the restriction of the
    multiplication above is a
    discrete symmetric space. The set $R:=T\times\{1,-1\}$ is a type
    $A_1$ 
    root system extended by the abelian group
    $G$ in the sense of \cite{yoshii-rootabelian} or \cite{myAbel}.


    Consider the two-element group
    $\calV:=\{1,-1\}$ with its action on $G$ characterized by 
    $-1g=-g$ for
    all $g\in G$. Set $\calA:=G\rtimes\calV$. Then $\calA$ acts on $T$
    via
    \begin{align*}
      (g,v).t=2g+vt.
    \end{align*}
    The map
    \begin{align*}
      \sbullet^\calA:~T\to\calA,~t\mapsto t^\calA=(t,-1)
    \end{align*}
    turns $\calA$ into a $T$-reflection group.

    In general, if $B$ is a group, $A$ is an abelian group, and 
    $f:B\times B\to A$ is a cocycle, then the set $A\times B$ with the
    multiplication given by
    \begin{align*}
      (a,b)(a',b')=\big(a+a'+f(b,b'),bb'\big)
    \end{align*}
    defines a group denoted by $A\times_f B$ which is a central
    extension of $B$. 

    The set $(G\wedge G)\times G\times\calV$ with the
    multiplication 
    \begin{align*}
      (l,g,v)(l',g',v')
      :=(l+l'+g\wedge(vg'),g+vg',vv')
    \end{align*}
    is a group. We denote it by 
    $(G\wedge G)\times_\wedge G\rtimes\calV$. 
    It can equally be interpreted as the semidirect product
    of the Heisenberg group $(G\wedge G)\times_\wedge G$ with $\calV$
    or a central extension of $\calA$ by $G\wedge G$ with cocycle
    \begin{math}
      f:~\calA^2\to G\wedge G,~\big((g,v),(g',v')\big)\mapsto g\wedge(vg').
    \end{math}

    Set 
    \begin{align*}
      \sbullet^\calW:~T\to(G\wedge G)\times_\wedge G\rtimes\calV,
      t\mapsto t^\calW=(0,t,-1).
    \end{align*}
    Let $\calW$ be the subgroup of 
    $(G\wedge G)\times_\wedge G\rtimes\calV$
    generated by $T^\calW$. Then $(\calW,\sbullet^\calW)$ is a 
    $T$-reflection group with the action of $\calW$ on $T$ induced by
    the action of $\calA$ on $T$.
    \begin{defn}
      The group $\calW$ is called \emph{the Weyl group of $R$}.
    \end{defn}
    This definition of the Weyl group coincides with
    the definition of Weyl groups given in \cite{myAbel} if $G$ is
    free abelian and the one given in \cite{EAWG} if $G$ is finitely
    generated free abelian.  

    \section{The abelian 2-group case}                               %

    In this section we investigate the case where $G$ is an
    elementary abelian 2-group. So we may think of $G$ as a vector
    space over the Galois field $\F_2$ with two elements. From
    (\ref{eq:mult}) it
    immediately follows that $T$ has the trivial multiplication.

    Denote by
    $G\otimes_\sym G$ the
    subgroup of $G\otimes G$ generated by the elements of the set
    \begin{math}
      \{v\otimes v~|~v\in G\}.
    \end{math}
    The group homomorphism
    \begin{align*}
      G\otimes G\to G\otimes_\sym G
      \text{~~characterized by~~}
      g\otimes h\mapsto g\otimes h-h\otimes g
    \end{align*}
    factors through $G\wedge G$ giving a group homomorphism
    \begin{align*}
      \pi:~G\wedge G\to G\otimes_\sym G
      \text{~~characterized by~~}
      u\wedge v\mapsto u\otimes v-v\otimes u.
    \end{align*}
    If $B$ is an ordered basis of $G$ then 
    $\{b_1\wedge b_2~|~b_1<b_2\in B\}$ is a basis of $G\wedge G$. Its
    image under $\pi$ is linearly independent, so $\pi$ is injective.

    Define the map
    \begin{align*}
      \phi:~
      (G\wedge G)\times_\wedge G
      \to
      G\otimes_\sym G,~
      (t,g)\mapsto(\pi(t)+g\otimes g).
    \end{align*}

    \begin{thm}\label{thm:grpIso}
      The map $\phi$ is a group isomorphism such that
      $\phi(0,g)=g\otimes g$ for all $g\in G$.
    \end{thm}

    \begin{proof}
      To see that $\phi$ is a group homomorphism let 
      $s,t\in G\wedge G$ and $g,h\in G$. Then
      \begin{align*}
        \phi\big((s,g)(t,h)\big)
        &=
        \phi(s+t+g\wedge h,g+h)
        \\&=
        \phi(s)+\pi(t)+g\otimes h+h\otimes g+(g+h)\otimes(g+h)
        \\&=
        \phi(s)+\pi(t)+g\otimes g+h\otimes h
        =
        \phi(s,g)+\phi(t,h).
      \end{align*}
      It is clear that $\phi$ is surjective, since it has a generating
      set in its image. 

      Since we are working with characteristic 2, the map
      \begin{equation*}
        G\to G\vee G,~ v\mapsto v\vee v
      \end{equation*}
      is an injective group homomorphism. 
      We denote by $G\vee_\sym G$ 
      the additive subgroup of $G\vee G$ generated by
      $\{g\vee g~|~g\in G\}$. So we have a group isomorphism
      \begin{align*}
        G\vee_\sym G\to G.
      \end{align*}
      Its composition with the quotient
      homomorphism 
      \begin{math}
        G\otimes_\sym G\to G\vee_\sym G
      \end{math}
      yields a homomorphism
      \begin{equation*}
        \sqrt{\sbullet}:~G\otimes_\sym G\to G
        ~~~~~
        \text{satisfying}
        ~~~~~
        \sqrt{g\otimes g}=g.
      \end{equation*}
      It vanishes on the image of $\pi$, since
      \begin{align*}
        \sqrt{\pi(g\wedge h)}
        &=
        \sqrt{g\otimes h-h\otimes g} 
        =
        \sqrt{g\otimes h+h\otimes g} 
        \\
        &=
        \sqrt{(g+h)\otimes(g+h)} 
        -\sqrt{g\otimes g}-\sqrt{h\otimes h}
        \\
        &=
        g+h-g-h
        =
        0
      \end{align*}
      for all $g$ and $h\in G$.
      
      To show that $\phi$ is injective, 
      let $(t,v)\in\ker(\phi)$, so $\pi(t)=v\otimes v$. Taking the
      square root on both sides yields $0=v$. We conclude $\pi(t)=0$.
      Since $\pi$ is injective we obtain $t=0$.
    \end{proof}

    In this section the action of $\calV$ on $G$ is trivial, so the
    reflection group $\calA$ is given by the direct product
    $\calA=G\times\calV$. The Weyl group $\calW$ is given as the
    subgroup of
    \begin{math}
      \big((G\wedge G)\times_\wedge G\big)\times\calV
    \end{math}
    generated by the image of
    \begin{align*}
      \sbullet^\calW:~
      T\to\big((G\wedge G)\times_\wedge G\big)\times\calV,~
      t\mapsto t^\calW=(0,t,-1).
    \end{align*}
    Due to the preceding theorem, the Weyl group can also be given as
    the
    subgroup of
    \begin{math}
      \big(G\otimes_\sym G)\times\calV
    \end{math}
    generated by the image of
    \begin{align}\label{eq:weylGroup}
      \sbullet^\calW:~
      T\to(G\otimes_\sym G)\times\calV,~
      t\mapsto t^\calW=(t\otimes t,-1).
    \end{align}
    
    Let $F:=F(T\setminus\{0\})$ be the free vector space on the set
    $T\setminus\{0\}$ with the embedding 
    $\iota:T\setminus\{0\}\to F$. 
    The initial reflection group is given by
    \begin{math}
      \calU=F\times\calV 
    \end{math}
    with the map
    \begin{align*}
      T\to\calU,~
      t\mapsto
      \begin{cases}
        (\iota(t),-1)&\text{if $t\neq 0$}\\
        (0,-1)&\text{if $t=0$.}
      \end{cases}
    \end{align*}

    \begin{defn}
      A subset $M$ of $G$ is called 
      \emph{2-dependent}, if the elements of the set
      $\{g\otimes g~|~g\in M\}$ are linearly dependent in
      $G\otimes G$.
      The set $M$ is called 
      \emph{2-independent} if it is not 2-dependent.
    \end{defn}
    \begin{exmp}\label{exmp:2Dep}
      \begin{enumerate}
  \item[a)]
        A linearly independent subset $M$ of $G$ is 2-independent, due
        to the homomorphism $\sqrt{\sbullet}$ used in the proof of
        Theorem~\ref{thm:grpIso}. 
  \item[b)]
        Set $G=(\F_2)^2$. Then the set $M$ of all nonzero vectors in $G$
        is 2-independent, since the matrices
        \begin{equation*}
          \begin{pmatrix}
            1&0\\
            0&0
          \end{pmatrix},
          \begin{pmatrix}
            0&0\\
            0&1
          \end{pmatrix},
          \begin{pmatrix}
            1&1\\
            1&1
          \end{pmatrix}
        \end{equation*}
        are linearly independent.
  \item[c)]
        Set $G=(\F_2)^n$. 
        Any subset $M$ of $G$ with cardinality $|M|>\frac{n(n+1)}2$
        is 2-dependent, since 
        $\dim_{\F_2}\big(G\otimes_\sym G\big)=\frac{n(n+1)}2$.  
      \end{enumerate}
    \end{exmp}

    \begin{thm}\label{thm:caseTwo}
      The reflection morphism $\calU\to\calW$ is injective if and only
      if the set $T\setminus\{0\}$ is 2-independent in $G$.
    \end{thm}

    \begin{proof} We will use the form of the Weyl group $\calW$ given
      in
      (\ref{eq:weylGroup}).
      Suppose $\calU\to\calW$ is not injective. Then there is a
      non-trivial element in its kernel. This element can be written
      as 
      $\big(\sum_{i=1}^n\iota(t_i),\sigma\big)\in 
      G\otimes_\sym G\times\calV$ for distinct 
      elements
      $t_1,~t_2,\dots,t_n\in T\setminus\{0\}$ and 
      $\sigma\in \calV$. It follows that
      $\sigma=1$ and 
      \begin{math}
        \sum_{i=1}^n t_i\otimes t_i=0.
      \end{math}
      So $t_1$, $t_2,\dots,t_n$ are 2-dependent. This implies that
      $T\setminus\{0\}$ is 2-dependent.

      Conversely, suppose $T\setminus\{0\}$ is 2-dependent, say 
      \begin{math}
        \sum_{i=1}^n t_i\otimes t_i=0
      \end{math}
      for distinct elements $t_1,~t_2,\dots,t_n\in T\setminus\{0\}$ and
      $n\ge1$. Then
      $\big(\sum_{i=1}^nt_i,0\big)$ is a nontrivial element in the kernel of 
      $\calU\to\calW$.  
    \end{proof}
    Denote the reflection morphism $\calU\to\calW$ above by
    $\phi$. Then Example~\ref{exmp:2Dep} yields
    \begin{cor}\label{cor:kerphim}
      \begin{enumerate}
  \item 
      The map $\phi$ is injective if
      $T\setminus\{0\}$ is a basis of $G$.
  \item 
      The map $\phi$ is not injective if 
      $|T\setminus\{0\}|>\frac{n(n+1)}2$, where $n=\dim(G)$. 
  \item 
      If $T=G$, then $\phi$ is an isomorphism if and only if
      $\dim(G)\le2$. 
      \end{enumerate}
    \end{cor}

    \section{The free abelian case}                                 %

    In this section let $G$ be a free abelian group. 
    We will reduce the situation to that of the former section. More
    details can be found in \cite{myAbel}~Section~2, in particular in
    Construction~2.10. 

    Let
    $\calU$ be the initial $T$-reflection group and let $\calW$ be the
    Weyl group. The abelianizations $\calU^\ab$ and $\calW^\ab$ are
    $T^\ab$-reflection groups, where $T^\ab$ is the image of $T$ under
    the quotient homomorphism $G\to G_2=G/2G$. This is a discrete
    symmetric space with the trivial multiplication.
    More precisely $\calU^\ab$
    is the initial $T^\ab$-reflection group and
    $\calW^\ab=\big(G_2\wedge G_2\big)\times_\wedge G_2\times\calV$ is the
    Weyl group for the discrete symmetric space $T^\ab$.

    The $T$-reflection morphism $\calU\to\calW$ yields a
    $T^\ab$-morphism $\calU^\ab\to\calW^\ab$ and there is a group
    homomorphism $\psi$ making the following diagram commute:
    \begin{align*}
      \xymatrix{
        \ker(\calU\to\calW)
        \ar[rr]
        \ar[d]_\psi
        &&
        \calU
        \ar[rr]
        \ar[d]
        &&
        \calW
        \ar[d]
        \\
        \ker(\calU^\ab\to\calW^\ab)
        \ar[rr]
        &&
        \calU^\ab
        \ar[rr]
        &&
        \calW^\ab.
      }
    \end{align*}

    According to \cite{myAbel}~Theorem~4.16 the map $\psi$ is an
    isomorphism. With Theorem~\ref{thm:caseTwo} we have obtained the
    main result of this article:

    \begin{thm}\label{thm:main}
      The $T$-reflection homomorphism $\calU\to\calW$ is an
      isomorphism if and only if 
      $T^\ab\setminus\{0\}$ is 2-independent in $G/2G$. 
    \end{thm}

    Corollary~\ref{cor:kerphim} gives more information 
    in some specific cases. In particular, it confirms the observation
    made in \cite{myWeyl} and \cite{azamA1} that $\calU\to\calW$ is
    not always injective.
    If $n$ is the rank of $G$ then testing for 2-dependence involves
    testing for linear dependence of $|T\setminus\{0\}|$ vectors in an
    $\frac{n(n+1)}2$-dimensional vector space over the Galois field
    $\F_2$. 
    This is more practical than testing for the existence of a
    so-called non-trivial integral collection
    according to \cite{azamA1}~Theorem~5.16. This theorem also
    requires $G$ to be finitely generated, a hypothesis that we don't
    require for our 
    Theorem~\ref{thm:main}.

    The hypotheses ``free'' for $G$ is only used to apply Theorem~4.16 of
    \cite{myAbel}. We would be interested in understanding if it could be
    weakened to ``torsion free'', ``involution free'' or even omitted
    completely.

     \bibliographystyle{alpha}
     \bibliography{../../references/references}
   \end{document}